\documentclass{article}
\usepackage{amsmath,epsfig}
\title{Correction to: ``Blind maximum likelihood separation of a linear-quadratic mixture''}
\author{Shahram Hosseini, Yannick Deville \vspace*{10pt} \\ \small
Laboratoire d'Astrophysique de Toulouse-Tarbes, Universit\'{e} de
Toulouse,  CNRS , \\ \small 14 avenue Edouard Belin - 31400
Toulouse - France. \\ \small shosseini@ast.obs-mip.fr ,
ydeville@ast.obs-mip.fr}
\date{}
\def\thesection{\Roman{section}}
\renewcommand{\theequation}{\thesection.\arabic{equation}}
\begin{document}
\maketitle
\begin{abstract}
An error occurred in the computation of a gradient in
\cite{hoss04}. The equations (20) in Appendix and (17) in the text
were not correct. The current paper presents the correct version
of these equations.
\end{abstract}
\section{Summary of \cite{hoss04}}
In \cite{hoss04} (see Appendix for an authors' version of this
article), we proposed a maximum likelihood approach for blindly
separating a linear-quadratic mixture defined by (Eq. (2) in
\cite{hoss04}):
\begin{equation}
x_1=s_1-l_1s_2-q_1s_1s_2 \;\;\;\;\;\;\;\;\;\;\;\;\;\;\;\;\;
x_2=s_2-l_2s_1-q_2s_1s_2 \label{mixture_model_e.eq}
\end{equation}
where $s_1$ and $s_2$ are two independent sources. The
log-likelihood for $N$ samples of the mixed signals $x_1$ and
$x_2$ reads (Eq. (12) in \cite{hoss04}):
\begin{equation}
L=E_t[\log{f_{S_1}(s_1(t))}]+E_t[\log{f_{S_2}(s_2(t))}]
-E_t[\log{|J(s_1(t),s_2(t))|}] \label{finalcost_e.eq}
\end{equation}
where $E_t[.]$ represents the time average operator on the $N$
samples, $f_{s_1}(.)$ and $f_{s_2}(.)$ are the probability density
functions (pdf) of the sources $s_1$ and $s_2$ and $J$ is the
Jacobian of the mixture which reads (Eq. (4) in \cite{hoss04})
\begin{equation}
J=1-l_1l_2-(q_2+l_2q_1)s_1-(q_1+l_1q_2)s_2. \label{jacobian_e.eq}
\end{equation}
Maximizing the log-likelihood requires that its gradient with
respect to the parameter vector ${\bf w}=[l_1,l_2,q_1,q_2]$, {\it
i.e.} $\frac{\partial L}{\partial {\bf w}}$, vanishes. Defining
the score functions of the two sources as (Eq. (13) in
\cite{hoss04})
\begin{displaymath}
\psi_{i}(u)=-\frac{\partial \log{f_{S_i}(u)}}{\partial u}
\;\;\;\;\; i=1,2
\end{displaymath}
we can write (Eq. (14) in \cite{hoss04})
\begin{equation}
\frac{\partial L}{\partial {\bf w}}=-E_t[\psi_1(s_1)\frac{\partial
s_1}{\partial {\bf w}}]-E_t[\psi_2(s_2)\frac{\partial
s_2}{\partial {\bf w}}]-E_t[\frac{1}{J}\frac{\partial J}{\partial
{\bf w}}] \label{gradient_e.eq}
\end{equation}
Rewriting (\ref{mixture_model_e.eq}) in the vector form ${\bf
x}={\bf f}({\bf s},{\bf w})$ and considering ${\bf w}$ as the
independent variable and ${\bf s}$ as the dependent variable, we
can write, using implicit differentiation (Eq. (15) in
\cite{hoss04})
\begin{equation}
{\bf 0}=\frac{\partial {\bf f}}{\partial{\bf s}}\frac{\partial
{\bf s}}{\partial{\bf w}}+\frac{\partial{\bf f}}{\partial{\bf w}}
\end{equation}
which yields (Eq. (16) in \cite{hoss04})
\begin{equation}
\frac{\partial {\bf s}}{\partial{\bf w}}=-(\frac{\partial {\bf
f}}{\partial{\bf s}})^{-1}\frac{\partial{\bf f}}{\partial{\bf w}}
\label{dsdw_e.eq}
\end{equation}
Note that $\frac{\partial {\bf f}}{\partial{\bf s}}$ is the
Jacobian matrix of the mixing model. Considering
(\ref{mixture_model_e.eq}), we can write (Appendix in
\cite{hoss04})
\\ $\frac{\partial {\bf f}}{\partial {\bf s}}= \left (
\begin{array}{cc}
1-q_1s_2 & -l_1-q_1s_1 \\ -l_2-q_2s_2 & 1-q_2s_1
\end{array}
\right )$ and $\frac{\partial {\bf f}}{\partial {\bf w}}= \left (
\begin{array}{cccc}
-s_2 & 0 & -s_1s_2 & 0 \\ 0 & -s_1 & 0 & -s_1s_2
\end{array}
\right )$, which implies, from (\ref{dsdw_e.eq})
\begin{equation}
\frac{\partial {\bf s}}{\partial {\bf w}}= \frac{-1}{J} \left (
\begin{array}{cc}
1-q_2s_1 & l_1+q_1s_1 \\ l_2+q_2s_2 & 1-q_1s_2
\end{array}
\right ).  \left (
\begin{array}{cccc}
-s_2 & 0 & -s_1s_2 & 0 \\ 0 & -s_1 & 0 & -s_1s_2
\end{array}
\right )  \label{dsdw_new_e.eq}
\end{equation}
and yields (Eq. (19) in \cite{hoss04})
\begin{eqnarray}
\frac{\partial s_1}{\partial{\bf
w}}=\frac{1}{J}\mbox{\huge$[$}(1-q_2s_1)s_2 \;,\; (l_1+q_1s_1)s_1
\;, (1-q_2s_1)s_1s_2 \;,\; (l_1+q_1s_1)s_1s_2 \mbox{\huge$]$}
\nonumber \\ \frac{\partial s_2}{\partial {\bf
w}}=\frac{1}{J}\mbox{\huge$[$}(l_2+q_2s_2)s_2 \;,\; (1-q_1s_2)s_1
\;, (l_2+q_2s_2)s_1s_2 \;,\; (1-q_1s_2)s_1s_2 \mbox{\huge$]$}
\label{ds1s2dw_e.eq}
\end{eqnarray}
Using (\ref{ds1s2dw_e.eq}), we obtain the first two terms of the
gradient (\ref{gradient_e.eq}). To obtain the third term, we need
to compute $\frac{\partial J}{\partial {\bf w}}$. This partial
derivative was computed inaccurately in \cite{hoss04} so that
Equations (20), and thus (17), in \cite{hoss04} are erroneous.
\section{Correct versions of Equations (20) and (17) in
\cite{hoss04}} \setcounter{equation}{0} In \cite{hoss04} we did
not consider the implicit dependence of $s_1$ and $s_2$ on ${\bf
w}$ and computed the derivative of $J$ with respect to ${\bf w}$
ignoring this dependence. Considering $J=g({\bf w},{\bf s}({\bf
w}))$, the correct equation for $\frac{\partial J}{\partial {\bf
w}}$ reads
\begin{equation}
\frac{\partial J}{\partial {\bf w}}=\frac{\partial J}{\partial
{\bf w}}_{\mid_{{\bf s}\; cte}}+\frac{\partial J}{\partial {\bf
s}}\frac{\partial {\bf s}}{\partial {\bf w}} \label{djdw_e.eq}
\end{equation}
Equation (20) in \cite{hoss04} only corresponded to the first term
on the right side of the above relation which reads, following
(\ref{jacobian_e.eq}), as:
\begin{equation}
\frac{\partial J}{\partial{\bf w}}_{\mid_{{\bf s}\;
cte}}=-\mbox{\huge$[$}l_2+q_2s_2 , l_1+q_1s_1 ,
 l_2s_1+s_2, s_1+l_1s_2 \mbox{\huge$]$}
\label{djdw_scte_e.eq}
\end{equation}
We now compute the gradient (\ref{djdw_e.eq}) entirely.
Considering (\ref{jacobian_e.eq}), we can write
\begin{equation}
\frac{\partial J}{\partial {\bf s}}=-[q2+l_2 q_1, q_1+l_1 q_2]
\label{djdws_e.eq}
\end{equation}
Using (\ref{djdw_e.eq}), (\ref{djdw_scte_e.eq}),
(\ref{djdws_e.eq}) and (\ref{dsdw_new_e.eq}) we finally obtain the
following equation which must replace the equation (20) in
\cite{hoss04}
\begin{eqnarray}
\frac{\partial J}{\partial {\bf w}} =[
-(l_2+q_2s_2)-(q_2+l_2q_1)(1-q_2s_1)s_2/J-(q_1+l_1q_2)(l_2+q_2s_2)s_2/J,
\nonumber
\\
-(l_1+q_1s_1)-(q_1+l_1q_2)(1-q_1s_2)s_1/J-(q_2+l_2q_1)(l_1+q_1s_1)s_1/J,
\nonumber \\
-(l_2s_1+s_2)-(q_2+l_2q_1)(1-q_2s_1)s_1s_2/J-(q_1+l_1q_2)(l_2+q_2s_2)s_1s_2/J,
\nonumber \\
-(l_1s_2+s_1)-(q_1+l_1q_2)(1-q_1s_2)s_1s_2/J-(q_2+l_2q_1)(l_1+q_1s_1)s_1s_2/J]
\nonumber \\ \label{dJdW_new_e.eq}
\end{eqnarray}
Inserting (\ref{ds1s2dw_e.eq}) and (\ref{dJdW_new_e.eq}) in
(\ref{gradient_e.eq}), we obtain the following expression for the
gradient which must replace Equation (17) in \cite{hoss04}
\begin{eqnarray}
\frac{\partial L}{\partial {\bf w}} &=& -E_t \mbox{\huge$[$}
\mbox{\huge$($}\psi_1(s_1)(1-q_2s_1)s_2+\psi_2(s_2)(l_2+q_2s_2)s_2
\nonumber \\
&-&(l_2+q_2s_2)-(q_2+l_2q_1)(1-q_2s_1)s_2/J-(q_1+l_1q_2)(l_2+q_2s_2)s_2/J
\mbox{\huge$)$}/J, \nonumber
\\
&\mbox{\huge$($}&\psi_1(s_1)(l_1+q_1s_1)s_1+\psi_2(s_2)(1-q_1s_2)s_1
\nonumber \\
&-&(l_1+q_1s_1)-(q_1+l_1q_2)(1-q_1s_2)s_1/J-(q_2+l_2q_1)(l_1+q_1s_1)s_1/J\mbox{\huge$)$}/J,
\nonumber \\
&\mbox{\huge$($}&\psi_1(s_1)(1-q_2s_1)s_1s_2+\psi_2(s_2)(l_2+q_2s_2)s_1s_2
\nonumber \\
&-&(l_2s_1+s_2)-(q_2+l_2q_1)(1-q_2s_1)s_1s_2/J-(q_1+l_1q_2)(l_2+q_2s_2)s_1s_2/J\mbox{\huge$)$}/J,
\nonumber \\
&\mbox{\huge$($}&\psi_1(s_1)(l_1+q_1s_1)s_1s_2+\psi_2(s_2)(1-q_1s_2)s_1s_2
\nonumber \\
&-&(l_1s_2+s_1)-(q_1+l_1q_2)(1-q_1s_2)s_1s_2/J-(q_2+l_2q_1)(l_1+q_1s_1)s_1s_2/J)\mbox{\huge$)$}/J\mbox{\huge$]$}
\nonumber
\\ \label{gradient_new_e.eq}
\end{eqnarray}

\newpage
\section*{APPENDIX: Authors' version of \cite{hoss04}}
\def\thesection{\arabic{section}}
\renewcommand{\theequation}{\arabic{equation}}
\setcounter{section}{0} \setcounter{equation}{0}
\section{Introduction}
\label{sec:intro}  It is well known that the independence
hypothesis is not sufficient for separating general nonlinear
mixtures because of the very large indeterminacies which make the
nonlinear BSS problem ill-posed. A natural idea for reducing the
indeterminacies is to constrain the structure of mixing and
separating models to belong to a certain set of transformations.
This supplementary constraint can be viewed as a regularization of
the initially ill-posed problem.

In this paper, we study a linear-quadratic mixture model which may
be considered as the simplest (nonlinear) version of a general
polynomial model. Our main aim is to develop an approach which can
be easily extended to higher-order polynomial models. Hence, we
propose a recurrent separating structure whose realization does
not require the knowledge of the explicit form of the inverse of
the mixing model. We develop a rigorous method to identify the
parameters of the separating structure in a maximum likelihood
framework. The algorithm is developed so that the inverse of the
mixing structure is not required to be known. Thus, it can be
extended to more general polynomial mixtures.
\section{mixing and separating models}
\label{sec:structure} Suppose $u_1$ and $u_2$ are two independent
random signals. Given the following nonlinear instantaneous
mixture model
\begin{equation}
x_i=a_{i1}u_1+a_{i2}u_2+b_{i}u_1u_2 \;\;\;\;\;i=1,2
\label{firstmodel.eq}
\end{equation}
we would like to estimate $u_1$ and $u_2$ up to a permutation and
a scaling factor (and possibly an additive constant). For
simplicity, let's denote $s_1=a_{11}u_1$ and $s_2=a_{22}u_2$.
$s_1$ and $s_2$ will be referred to as the {\it sources} in the
following. (\ref{firstmodel.eq}) can be rewritten as
\begin{eqnarray}
x_1=s_1-l_1s_2-q_1s_1s_2 \nonumber \\ x_2=s_2-l_2s_1-q_2s_1s_2
\label{mixture_model.eq}
\end{eqnarray}
in which $l_1=-a_{12}/a_{22}$ and $l_2=-a_{21}/a_{11}$ represent
the linear contributions of the sources in the mixture, and
$q_1=-b_1/(a_{11}a_{22})$ and $q_2=-b_2/(a_{11}a_{22})$ represent
the quadratic contributions. The negative signs are chosen for
simplifying the notations of the separating structure.

Solving the model (\ref{mixture_model.eq}) for $s_1$ and $s_2$
leads to the following two pairs of solutions, which may be
considered as two direct separating structures:
\begin{eqnarray}
({\cal s}_1, {\cal s}_2)_1=((-b_1+\sqrt{\Delta_1})/2a_1,
(-b_2+\sqrt{\Delta_2})/2a_2) \nonumber \\ ({\cal s}_1, {\cal
s}_2)_2=((-b_1-\sqrt{\Delta_1})/2a_1, (-b_2-\sqrt{\Delta_2})/2a_2)
\label{inverse.eq}
\end{eqnarray}
where $\Delta_i=b_i^2-4a_ic_i$, $a_1=q_2+l_2q_1$,
$a_2=q_1+l_1q_2$, $b_1=q_1x_2-q_2x_1+l_1l_2-1$,
$b_2=q_2x_1-q_1x_2+l_1l_2-1$, $c_1=x_1+l_1x_2$ and
$c_2=x_2+l_2x_1$. It can be easily verified that
$\Delta_1=\Delta_2=J^2$, where $J$ is the Jacobian of the mixing
model (\ref{mixture_model.eq}) and reads
\begin{equation}
J=1-l_1l_2-(q_2+l_2q_1)s_1-(q_1+l_1q_2)s_2 \label{jacobian.eq}
\end{equation}
According to the variation domain of the two sources, three
different cases may be considered:

1) $J<0$ for all the values of $s_1$ and $s_2$. In this case
(\ref{inverse.eq}) becomes:
\begin{equation}
({\cal s}_1, {\cal s}_2)_1 = (s_1,s_2) \label{pair1.eq}
\end{equation}
\begin{equation}
 ({\cal s}_1, {\cal
s}_2)_2 =
(-\frac{q_1+l_1q_2}{q_2+l_2q_1}s_2-\frac{l_1l_2-1}{q_2+l_2q_1},
-\frac{q_2+l_2q_1}{q_1+l_1q_2}s_1-\frac{l_1l_2-1}{q_1+l_1q_2})
\label{pair2.eq}
\end{equation}
Thus, the first direct separating structure in (\ref{inverse.eq})
leads to the actual sources and the second direct separating
structure leads to another solution, equivalent to the first one
up to a permutation, a scaling factor, and an additive constant.

\begin{figure}
\centerline{\psfig{figure=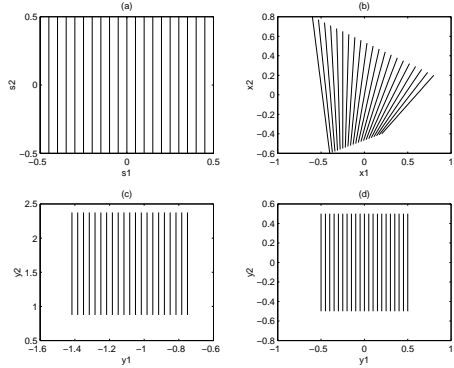,width=60mm}}
\caption{\label{case2.fig} Case when $J>0$ for all the source
values. Distribution of (a) sources, (b) mixtures, (c) output of
the first direct separating structure, (d) output of the second
direct separating structure.}
\end{figure}

2) $J>0$ for all the values of $s_1$ and $s_2$. In this case, the
first structure leads to the permuting solution, defined by
(\ref{pair2.eq}), and the second structure to the actual sources
$(s_1, s_2)$. An example is shown in Fig. \ref{case2.fig} for the
numerical values $l_1=-0.2$, $l_2=0.2$, $q_1=-0.8$, $q_2=0.8$ and
$s_i \in [-0.5, 0.5]$.

\begin{figure}
\centerline{\psfig{figure=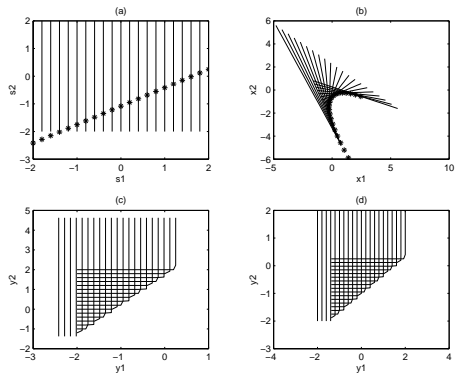,width=60mm}}
\caption{\label{case3.fig} Case when $J>0$ for some values of the
sources and \mbox{$J<0$} for the other values. Distribution of (a)
sources, (b) mixtures, (c) output of the first direct separating
structure, (d) output of the second direct separating structure.}
\end{figure}

3) $J>0$ for some values of the sources and $J<0$ for the other
values. In this case, each structure leads to the non-permuted
sources (\ref{pair1.eq}) for some values of the observations and
to the permuted sources (\ref{pair2.eq}) for the other values. An
example is shown in Fig. \ref{case3.fig} (with the same
coefficients as in the second case, but for \mbox{$s_i \in [-2,
2]$}). The permutation effect is clearly visible in the figure.
One may also remark that the straight line $J=0$ in the source
plane is mapped to a conic section in the observation plane (shown
by asterisks).

Thus, it is clear that the direct structures may be used for
separating the sources if the Jacobian of the mixing model is
always negative or always positive, {\it i.e.} for all the source
values. Otherwise, although the sources are separated {\it sample
by sample}, each retrieved signal contains samples of the two
sources. This problem arises because the mixing model
(\ref{mixture_model.eq}) is not bijective. This theoretically
insoluble problem should not discourage us. In fact, our final
objective is to extend the idea developed in the current study to
more general polynomial models which will be used to approximate
the nonlinear mixtures encountered in the real world. If these
real-world nonlinear models are bijective, we can logically
suppose that the coefficients of their polynomial approximations
take values which make them bijective on the variation domains of
the sources. Thus, in the following, we suppose that the sources
and the mixture coefficients have numerical values ensuring that
the Jacobian $J$ of the mixing model has a constant sign.

\begin{figure}
\centerline{\psfig{figure=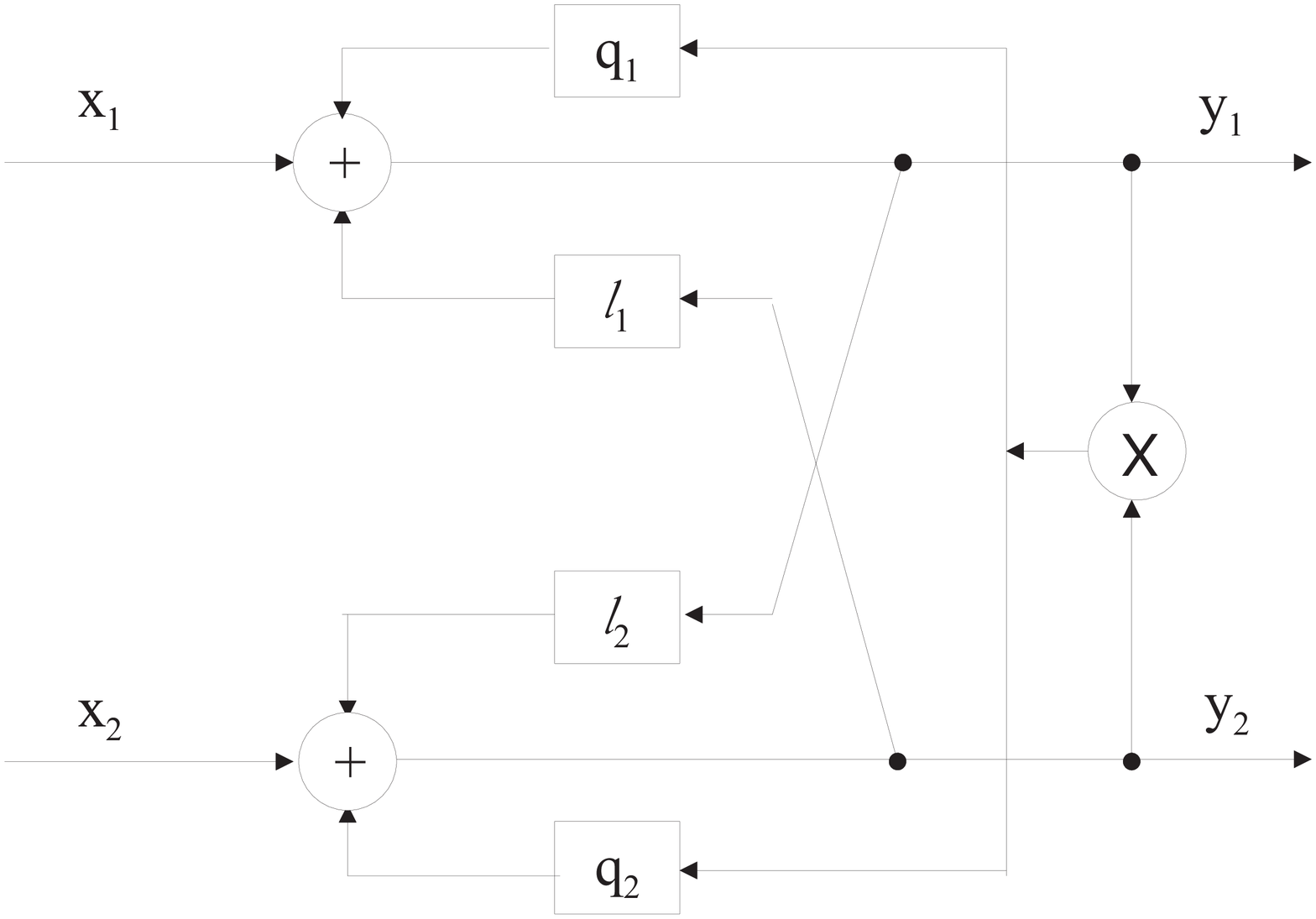,width=50mm}}
\caption{\label{model.fig} Recurrent separating structure.}
\end{figure}

The natural idea to separate the sources is to form a direct
separating structure using any of the equations in
(\ref{inverse.eq}), and to identify the parameters $l_1$, $l_2$,
$q_1$ and $q_2$ by optimizing an independence measuring criterion.
Although this approach may be used for our special mixing model
(\ref{mixture_model.eq}), as soon as a more complicated polynomial
model is considered, the solutions $({\cal s}_1, {\cal s}_2)$ can
no longer be determined so that the generalization of the method
to arbitrary polynomial models seems impossible. To avoid this
limitation, we propose a recurrent structure shown in Fig.
\ref{model.fig}. Note that, for $q_1=q_2=0$, this structure is
reduced to the basic H\'{e}rault-Jutten network. It may be checked
easily that, for fixed observations defined by
(\ref{mixture_model.eq}), $y_1=s_1$ and $y_2=s_2$ corresponds to a
steady state for the structure in Figure \ref{model.fig}.

The use of this recurrent structure is more promising because it
can be easily generalized to arbitrary polynomial models. However,
the main problem with this structure is its stability. In fact,
even if the mixing model coefficients are exactly known, the
computation of the structure outputs requires the realization of
the following recurrent iterative model
\begin{eqnarray}
y_1(n+1)=x_1+l_1y_2(n)+q_1y_1(n)y_2(n) \nonumber \\
y_2(n+1)=x_2+l_2y_1(n)+q_2y_1(n)y_2(n) \label{recurrent.eq}
\end{eqnarray}
where a loop on $n$ is performed for each couple of observations
$(x_1,x_2)$ until convergence is achieved.

It can be shown that this model is locally stable at the
separating point $(y_1,y_2)=(s_1,s_2)$, if and only if the
absolute values of the two eigenvalues of the Jacobian matrix of
(\ref{recurrent.eq}) are smaller than one. In the following, we
suppose that this condition is satisfied.
\section{Maximum likelihood estimation of the model parameters}
\label{sec:ML} Let $f_{S_1,S_2}(s_1,s_2)$ be the joint pdf of the
sources, and assume that the mixing model is bijective so that the
Jacobian of the mixing model has a constant sign on the variation
domain of the sources. The joint pdf of the observations can be
written as
\begin{equation}
f_{X_1,X_2}(x_1,x_2)=\frac{f_{S_1,S_2}(s_1,s_2)}{|J(s_1,s_2)|}
\label{pdf1.eq}
\end{equation}
Taking the logarithm of (\ref{pdf1.eq}), and considering the
independence of the sources, we can write:
\begin{equation}
\log{f_{X_1,X_2}(x_1,x_2)}=\log{f_{S_1}(s_1)}+ \log{f_{S_2}(s_2)}
-\log{|J(s_1,s_2)|} \label{logf1.eq}
\end{equation}
Given N samples of the mixtures $X_1$ and $X_2$, we want to find
the maximum likelihood estimator for the mixture parameters ${\bf
w}=[l_1,l_2,q_1,q_2]$. This estimator is obtained by maximizing
the joint pdf of all the observations (supposing that the
parameters in ${\bf w}$ are constant), which is equal to
\begin{equation}
E=f_{X_1,X_2}(x_1(1),x_2(1),\cdots,x_1(N),x_2(N))
\end{equation}
If $s_1(t)$ and $s_2(t)$ are two i.i.d. sequences, $x_1(t)$ and
$x_2(t)$ are also i.i.d. so that
$E=\prod_{i=1}^N{f_{X_1,X_2}(x_1(i),x_2(i))}$ and
$\log{E}=\sum_{i=1}^N{\log{f_{X_1,X_2}(x_1(i),x_2(i))}}$. The cost
function to be maximized can be defined as $L=\frac{1}{N}\log{E}$,
which will be denoted using the temporal averaging operator
$E_t[.]$ as
\begin{equation}
L= E_t[\log{f_{X_1,X_2}(x_1(t),x_2(t))}]
\end{equation}
Using (\ref{logf1.eq}):
\begin{equation}
L=E_t[\log{f_{S_1}(s_1(t))}]+E_t[\log{f_{S_2}(s_2(t))}]
-E_t[\log{|J(s_1(t),s_2(t))|}] \label{finalcost.eq}
\end{equation}
Maximizing this cost function requires that its gradient with
respect to the parameter vector ${\bf w}$, {\it i.e.}
$\frac{\partial L}{\partial {\bf w}}$, vanishes. Defining the
score functions of the two sources as
\begin{equation}
\psi_{i}(u)=-\frac{\partial \log{f_{S_i}(u)}}{\partial u}
\;\;\;\;\; i=1,2
\end{equation}
and considering that $\frac{\partial \log{|J|}}{\partial {\bf w}}=
\frac{1}{J}\frac{\partial J}{\partial {\bf w}}$, we can write
\begin{equation}
\frac{\partial L}{\partial {\bf w}}=-E_t[\psi_1(s_1)\frac{\partial
s_1}{\partial {\bf w}}]-E_t[\psi_2(s_2)\frac{\partial
s_2}{\partial {\bf w}}]-E_t[\frac{1}{J}\frac{\partial J}{\partial
{\bf w}}] \label{gradient.eq}
\end{equation}
Rewriting (\ref{mixture_model.eq}) in the vector form ${\bf
x}={\bf f}({\bf s},{\bf w})$ and considering ${\bf w}$ as the
independent variable and ${\bf s}$ as the dependent variable, we
can write, using implicit differentiation
\begin{equation}
{\bf 0}=\frac{\partial {\bf f}}{\partial{\bf s}}\frac{\partial
{\bf s}}{\partial{\bf w}}+\frac{\partial{\bf f}}{\partial{\bf w}}
\end{equation}
which yields
\begin{equation}
\frac{\partial {\bf s}}{\partial{\bf w}}=-(\frac{\partial {\bf
f}}{\partial{\bf s}})^{-1}\frac{\partial{\bf f}}{\partial{\bf w}}
\label{dsdw.eq}
\end{equation}
Note that $\frac{\partial {\bf f}}{\partial{\bf s}}$ is the
Jacobian matrix of the mixing model. Using (\ref{gradient.eq}) and
(\ref{dsdw.eq}), the gradient of the cost function $L$ with
respect to the parameter vector ${\bf w}$ is equal to (see the
appendix for the computation details)
\begin{eqnarray}
\frac{\partial L}{\partial {\bf w}} = -E_t \mbox{\huge$[$}
\mbox{\huge$($}\psi_1(s_1)(1-q_2s_1)s_2+\psi_2(s_2)(l_2+q_2s_2)s_2
-(l_2+q_2s_2) \mbox{\Large$)$}/J, \nonumber \\
\mbox{\Large$($}\psi_1(s_1)(l_1+q_1s_1)s_1+\psi_2(s_2)(1-q_1s_2)s_1
-(l_1+q_1s_1)\mbox{\Large$)$}/J, \nonumber \\
\mbox{\Large$($}\psi_1(s_1)(1-q_2s_1)s_1s_2+\psi_2(s_2)(l_2+q_2s_2)s_1s_2
-(l_2s_1+s_2)\mbox{\Large$)$}/J, \nonumber \\
\mbox{\Large$($}\psi_1(s_1)(l_1+q_1s_1)s_1s_2+\psi_2(s_2)(1-q_1s_2)s_1s_2
-(s_1+l_1s_2)\mbox{\Large$)$}/J\mbox{\huge$]$} \label{dLdw.eq}
\end{eqnarray}
In practice, the actual sources and their density functions are
unknown and will be replaced by the reconstructed sources, {\it
i.e.} by the outputs of the separating structure of Fig
\ref{model.fig}, $y_i$, in an iterative algorithm. The score
functions of the reconstructed sources can be estimated by any of
the existing parametric or non-parametric methods. In our work, we
used a kernel estimator based on third-order cardinal splines.
Using (\ref{dLdw.eq}), the cost function (\ref{finalcost.eq}) can
be maximized by a gradient ascent algorithm which updates the
parameters by the rule ${\bf w}(n+1)={\bf w}(n)+\mu \frac{\partial
L}{\partial {\bf w}}$. The learning rate parameter $\mu$ must be
chosen carefully to avoid the divergence of the algorithm. Note
that the algorithm does not require the knowledge of the explicit
inverse of the mixing model (direct separating structures
(\ref{inverse.eq})). Hence, it can be easily extended to more
general polynomial mixing models.
\section*{\center{Appendix: details of gradient
computation}} 
 Considering (\ref{mixture_model.eq}), we can write
\\ $\frac{\partial {\bf f}}{\partial {\bf s}}= \left (
\begin{array}{cc}
1-q_1s_2 & -l_1-q_1s_1 \\ -l_2-q_2s_2 & 1-q_2s_1
\end{array}
\right )$ and $\frac{\partial {\bf f}}{\partial {\bf w}}= \left (
\begin{array}{cccc}
-s_2 & 0 & -s_1s_2 & 0 \\ 0 & -s_1 & 0 & -s_1s_2
\end{array}
\right )$, which implies, from (\ref{dsdw.eq})
\begin{equation}
\frac{\partial {\bf s}}{\partial {\bf w}}= \frac{-1}{J} \left (
\begin{array}{cc}
1-q_2s_1 & l_1+q_1s_1 \\ l_2+q_2s_2 & 1-q_1s_2
\end{array}
\right ).  \left (
\begin{array}{cccc}
-s_2 & 0 & -s_1s_2 & 0 \\ 0 & -s_1 & 0 & -s_1s_2
\end{array}
\right ) \nonumber
\end{equation}
which yields
\begin{eqnarray}
\frac{\partial s_1}{\partial{\bf
w}}=\frac{1}{J}\mbox{\huge$[$}(1-q_2s_1)s_2 \;,\; (l_1+q_1s_1)s_1
\;, (1-q_2s_1)s_1s_2 \;,\; (l_1+q_1s_1)s_1s_2 \mbox{\huge$]$}
\nonumber \\ \frac{\partial s_2}{\partial {\bf
w}}=\frac{1}{J}\mbox{\huge$[$}(l_2+q_2s_2)s_2 \;,\; (1-q_1s_2)s_1
\;, (l_2+q_2s_2)s_1s_2 \;,\; (1-q_1s_2)s_1s_2 \mbox{\huge$]$}
\label{ds1s2dw.eq}
\end{eqnarray}
Considering (\ref{jacobian.eq})
\begin{eqnarray}
\frac{\partial J}{\partial{\bf w}}=-\mbox{\huge$[$}l_2+q_2s_2 ,
l_1+q_1s_1 ,
 l_2s_1+s_2, s_1+l_1s_2 \mbox{\huge$]$}
\label{djdw.eq}
\end{eqnarray}
(\ref{dLdw.eq}) follows directly from (\ref{gradient.eq}),
(\ref{ds1s2dw.eq}) and (\ref{djdw.eq}).

\begin{thebibliography}{9}
\bibitem{hoss04} S. Hosseini, Y. Deville, Blind maximum likelihood separation of a linear-quadratic mixture,
Proceedings of the Fifth International Conference on Independent
Component Analysis and Blind Signal Separation (ICA 2004), pp.
694-701, ISSN 0302-9743, ISBN 3-540-23056-4, Springer-Verlag, vol.
LNCS 3195, Granada, Spain, Sept. 22-24, 2004.
\end{thebibliography}
\end{document}